\documentclass[12pt]{amsart}
\usepackage{amsmath, url, hyperref, amssymb}
\usepackage{graphicx}
\usepackage{caption}
\usepackage{subcaption}
\usepackage{booktabs}

\newtheorem{theorem}{Theorem}[section] 

\newtheorem{conjecture}[theorem]{Conjecture}

\newtheorem{example}{Example}

\newcommand{\FF}{\mathbb{F}}

\newcommand{\QQ}{\mathbb{Q}}

\newcommand{\ZZ}{\mathbb{Z}}
\newcommand{\frakp}{\mathfrak{p}}
\newcommand{\eqdef}{{:=}}

\def\ovl{\overline}
\def\ra{\rightarrow}

\DeclareMathOperator{\NS}{NS}
\DeclareMathOperator{\Br}{Br}
\DeclareMathOperator{\disc}{disc}

\newlength{\graphwidthtwo}
\newlength{\graphwidthone}
\title[N\'eron-Severi ranks of K3 surfaces]{Variation of N\'eron-Severi ranks of reductions of K3 surfaces}

\date{\today}

\author{Edgar Costa}
\address{Courant Institute, NYU, 251 Mercer St.  New York, NY 10012, USA}
\email{edgarcosta@nyu.edu}

\author{Yuri Tschinkel}
\address{Courant Institute\\
                New York University \\
                New York, NY 10012 \\
                USA }
\email{tschinkel@cims.nyu.edu}

\address{Simons Foundation\\
160 Fifth Avenue\\
New York, NY 10010\\
USA}

\begin{document}
\maketitle

\begin{abstract}
We study the behavior of geometric Picard ranks of K3 surfaces over $\mathbb Q$ 
under reduction modulo primes. 
We compute these ranks for reductions of smooth quartic surfaces modulo all primes $p<2^{16}$ 
in several representative examples and investigate the resulting statistics.
\end{abstract}

\maketitle

\section{Introduction}
\label{sect:intro}

Let $k$ be a number field and $X$ a K3 surface over $k$, i.e., a geometrically smooth projective simply-connected 
surface with trivial canonical class, for example, a smooth quartic hypersurface in $\mathbb P^3$.  
Let $\frakp$ be a finite place of $k$ where $X$ has good reduction $X_\frakp$.
Let $\overline{X}$ (resp. $\overline{X}_\frakp$) be
the base change of $X$ (respectively, $X_\frakp$) 
to the algebraic closure of $k$ (respectively, 
of the residue field of $\frakp$), 
and let $\rho(\ovl{X})$ and $\rho(\ovl{X}_\frakp)$ be the 
ranks of the corresponding N\'eron-Severi groups  $\NS(\ovl{X})$ and $\NS (\ovl{X}_\frakp )$, i.e., the
geometric Picard ranks. There is a natural specialization homomorphism 
\begin{equation}
\label{eqn:special}
    s_{\mathfrak{p}} : \NS ( \ovl{X} ) \rightarrow \NS (\ovl{X}_\frakp ),
\end{equation}
which is injective (see, e.g., \cite[Proposition 6.2]{vluijktriangles}), thus 
$$
\rho(\ovl{X}) \leq \rho(\ovl{X}_\frakp).
$$
In fact, for all $\frakp$ of good reduction we have
\begin{equation}
\label{eqn:ineq}
\rho(\ovl{X}) + \eta(\ovl{X}) \le \rho(\ovl{X}_\frakp), 
\end{equation}
for some $\eta(\ovl{X}) \geq 0$, 
defined by \eqref{eqn:eta}. 
It is known that there exist infinitely many $\frakp$ such that equality occurs in (\ref{eqn:ineq});
furthermore, over some finite extension of $k$, 
the set of such primes has density one \cite[Theorem 1]{picardk3}.  
However, very little is known about the set of primes 
$$
\Pi_{\rm jump}(X):= \{ \frakp  :  \rho(\ovl{X}) + \eta(\ovl{X}) < \rho(\ovl{X}_\frakp) \},
$$ 
where the inequality \eqref{eqn:ineq} is strict. 

Information about $\Pi_{\rm jump}(X)$ can be converted into geometric statements:
{\em if} this set contains infinitely many primes of non-supersingular reduction, 
for all K3 surfaces over number fields with $\rho(X)=2,4$, 
then {\em all} K3 surfaces over algebraically closed fields 
of characteristic zero have infinitely many rational curves, 
by \cite{bht} and \cite{li-lie}. 

There are cases where 
$\Pi_{\rm jump}(X)$ is known to be infinite. For example, assume that $X$ is a Kummer surface, i.e., 
the resolution of singularities of the quotient $A/\iota$, 
where $A$ is an abelian surface, and 
$\iota : A\ra A$ the standard involution $\iota(a)=-a$. Then 
$$
\rho(\ovl{X})=\rho(\ovl{A})+16.
$$
Now assume that $A \sim C_1\times C_2$, i.e., is 
isogenous to a product of two elliptic curves.
Then 
\begin{itemize}
\item[(i)] $\rho(\ovl{X})\ge 18$, 
\item[(ii)] $\rho(\ovl{X})\ge 19$, if $C_1\sim C_2$, and 
\item[(iii)] $\rho(\ovl{X}) =  20$, if in addition, $C_1$ has complex multiplication by 
$E:=\mathbb Q(\sqrt{-d})$.    
\end{itemize}
In these extreme cases, the primes in $\Pi_{\rm jump}(X)$ can be understood as follows:
\begin{itemize}
\item if $\rho(\ovl{X}) \ge 19$, then $\frakp \in \Pi_{\rm jump}(X)$ provided $\frakp$ 
is a supersingular prime for $C_1$ (and thus $C_2$).
\end{itemize}

By a theorem 
of Elkies, there are infinitely many such primes \cite{elkies}, at least for elliptic curves over $\QQ$. 

In case (i), $\frakp \in \Pi_{\rm jump}(X)$ provided the reductions of 
$C_1$ and $C_2$ modulo $\frakp$ are isogenous. There are infinitely many such $\frakp$, by a recent theorem of Charles \cite{charles-inf}.

This motivates us to consider the asymptotic behavior of the proportion of primes in  
$\Pi_{\rm jump}(X)$:
\begin{equation}
\gamma(X,B) \eqdef 
\frac{\#\left\{ p \leq B :  \, p\in \Pi_{\rm jump}(X) \right\} }{\#\left\{p \leq B \right\}}.
\end{equation}

Returning to Kummer surfaces of the form $X\sim C\times C/\iota$, when the elliptic curve
$C$ does not have complex multiplication, so that $\rho(\ovl{X})=19$, 
the Lang-Trotter conjecture \cite{LT}, 
implies 
$$
\gamma(X,B) \sim \frac{c}{\sqrt{B}}, \quad B\ra \infty,
$$
for some constant $c > 0$.
The Lang-Trotter conjecture has attracted the attention of many experts; 
for a sample of results we
refer to \cite{elkiesLT,murtyLT,chantalfrancescoLT,baierLT} and to 
\cite{katz}, in the function field case.
If $C$ does have complex multiplication, then 
\begin{equation}
\label{eqn:12}
\gamma(X,B) \sim \frac{1}{2}, \quad B\ra \infty.
\end{equation}

Elsenhans and Jahnel conducted an extensive numerical investigation of 
Kummer surfaces over $\mathbb Q$, in particular of those with $\rho(\ovl{X})=17$ \cite{kummersurfaces}.
They computed $\rho(\overline{X}_p)$, for $p < 1000$,  for a large sample of surfaces $X$ with 
$\rho(\ovl{X})=17$, and observed that the proportion of such $X$
with $\rho(\ovl{X}_p) > 18$ is roughly $2/\sqrt{p}$.
In another direction, for some of the examples with $\rho(\overline{X}) = 18$, 
they were able to show that the density of $\Pi_{\rm jump}(X)$ is at least $1/2$.
The precise shape of asymptotic formulas for $\gamma(X,B)$  
for general Kummer surfaces $X\sim A/\iota$ 
is likely to depend on the Sato-Tate group $\mathrm{ST}_A$ 
of the abelian surface $A$, investigated in \cite{sutherland}.

More generally, the Kuga-Satake construction (see \cite{deligne}) relates a K3 surface $X$ to an 
abelian variety $A=A_X$ of dimension $2^{19}$. 
Knowing this abelian variety explicitly, 
in particular, its Picard group and its endomorphisms, 
would allow us to compute the Picard group of $X$, see \cite[Proposition 19]{hkt}. 
The jumping behavior of Picard ranks of K3 surfaces is therefore 
related to the jumping behavior 
on these abelian varieties, similar to the Kummer case above, thus 
should be controlled by a version of the Lang-Trotter conjecture. 
However, the Kuga-Satake construction is transcendental, and
we do not yet have sufficiently effective control over $A$, even over its field of definition, 
except in degree two \cite[Remark 9]{hkt}. 

In this note we report on a numerical study of the variation of Picard ranks of 
quartic K3 surfaces over $\QQ$,  with {\em small} $\rho(\ovl{X})$. 
For several representative examples, we compute 
$\rho(\overline{X})$ and $\rho(\overline{X}_p)$, for all $2 <p < 2^{16}$,
where $X$ has good reduction, and we calculate $\gamma(X,B)$, for $B < 2^{16}$.

We observe two different trends.
In examples where $\rho(\overline{X}) = 1$ and $\eta(\overline{X}) = 1$ 
we find evidence that 
$$
\gamma(X,B) \sim \frac{c_X}{\sqrt{B}}, \quad B\ra \infty,
$$ 
for some constant $c_X>0$. 
In other words, a prime $p$ is in $\Pi_{\rm jump}(X)$ with probability proportional to $1/\sqrt{p}$.
On our other examples, when $\rho(\overline{X}) = 2$ (and $\eta(\overline{X}) = 0$), 
we are lead to believe that 
$$
\liminf_{B\rightarrow \infty} \gamma(X,B) \geq \frac{1}{2},
$$ 
i.e., the primes at which the geometric Picard number jumps have density $\geq 1/2$. 
Our data strongly suggests that we are {\em not} in the same situation as in \eqref{eqn:12}, the plots in 
Figure~\ref{fig:picardnumbertwo} and Figure~\ref{fig:picardnumbertwomavg} are not consistent with
statistics for the splitting behavior of primes in quadratic extensions of 
$\QQ$. 

\

\noindent
{\bf Acknowledgments.} 
We are grateful to David Harvey, Brendan Hassett, Kiran Kedlaya, and Barry Mazur 
for useful discussions and collaboration on related questions. 
The first author was partially supported by FCT doctoral grant SFRH/BD/69914/2010.
The second author was supported by National Science Foundation
grants 0968318 and 1160859.

\section{Computing the Picard number of a K3 surface}
\label{sec:picard}

In this section, we explain our approach to the computation of Picard numbers of quartic K3 surfaces. 
Over a {\em finite field}, one only needs to compute the Hasse-Weil zeta function; which may be 
computationally expensive, but is achievable in bounded time. Over a {\em number field},  
computing the Picard number of an algebraic surface is a hard problem. 
For K3 surfaces, an effective version of the Kuga-Satake construction as in \cite{hkt} 
yields a theoretical algorithm, 
with {\em a priori} bounded running time, at least for degree-two K3 surfaces. 
In \cite[Section 8.6.]{poonen}  
the authors provide an alternative algorithm; another
algorithm, conditional on the Hodge conjecture for $X\times X$, is presented in   
\cite{picardk3}; these algorithms do not have {\em a priori} bounded running times.

In practice, one starts by establishing lower and upper bounds for $\rho(\ovl{X})$.
Lower bounds can be produced by exhibiting independent divisors on $\ovl{X}$,  
and upper bounds can be obtained via specialization to finite fields as in \eqref{eqn:special}.
This approach does not guarantee an answer in every case, but sometimes the bounds agree.  
In some cases, one can improve the upper bound 
by a careful analysis of the  specialization map.  
For example, if the lattice structure disagrees over two different 
specializations, or if some divisor class on $\ovl{X}_\frakp$ is not liftable, 
then the specializations cannot be surjective.  
This approach has its limitations, as one cannot in general 
expect that there exist places $\frakp$ 
such that $\rho(\ovl{X}_\frakp) \leq \rho(\ovl{X}) + 1$. 
An overview of these techniques can be found in \cite[Chapter 7]{2lectk3}.

In \cite{picardk3}, Charles proved a general theorem about the jumping behavior of 
Picard ranks under specialization: Let $E_X$ be the endomorphism algebra of the 
Hodge structure underlying the transcendental lattice $T_X$ of $X$; 
it is known that $E_X$ is a field, which is either totally real or a CM-field (see, e.g., \cite{zarhin}).
In the latter case, one says that $X$ has complex multiplication.
By \cite[Theorem 1]{picardk3}, there are two possibilities,
\begin{equation}
  \label{finiteplacepicard}
  \rho(\ovl{X}_\frakp) \geq 
  \begin{cases}
      \rho(\ovl{X}) &\text{if $E_X$ is a CM-field or $\dim_{E_X}(T_X)$ is even},\\
      \rho(\ovl{X}) + [E_X:\QQ] &\text{if $E_X$ is totally real field and 
                                      $\dim_{E_X}(T_X)$ is odd}. 
  \end{cases}
\end{equation}
We define 
\begin{equation}
\label{eqn:eta}
\eta(\ovl{X}) := 0 \quad \text{ or } \quad [E_X:\QQ],
\end{equation}
depending on which case we are in. 


We turn to finite fields. Let $X$ be a smooth projective surface over $\FF_q$.
The Weil conjectures tell us that the Hasse-Weil zeta function has the form
\begin{equation}
  \label{zeta}
  Z(X, T) \eqdef \exp\left( \sum_{m=1} ^{\infty} \frac{ \# X(\FF_{q^m})}{m} t^m \right) 
= \frac{P_1 (X, t) P_3(X, t)}{(1-t) P_2(X, t) (1-q^2 t)},
\end{equation}
where 
$$
P_i(X, t):=\det\left(1 - t \, \mathrm{Fr}_i | H_{\text{et}} ^i (\overline{X},\QQ_{\ell})\right)\in \mathbb Z[t]
$$
have reciprocal roots of absolute value $q^{i/2}$, and $\mathrm{Fr}_i$ are the Frobenius automorphisms.
The Artin-Tate conjecture relates the N\'eron-Severi group of $X$ with $P_2(X,t)$:
\begin{conjecture}

\

\begin{itemize}
\item (Tate Conjecture)
$\rho(X)$ equals the multiplicity of $q$ as a reciprocal root of $P_2(X, t)$.
\item (Artin-Tate Conjecture)
Let $\Br(X)$ be the Brauer group of $X$  and  
$$
\alpha(X) := \chi(X,\mathcal{O}_X) - 1 + \dim(\operatorname{Pic} ^0 (X)).
$$
Then 
  \begin{equation*}
    \lim_{s \rightarrow 1} \frac{ P_2(X,q^{-s})}{(1-q^{1-s})^{\rho(X)} } = \frac{(-1)^{\rho(X) -1} \# \Br(X)\cdot \disc( \NS(X))}{q^{\alpha(X)} ( \# \NS(X)_{\text{tors}})^2}.
  \end{equation*}
\end{itemize}
\end{conjecture}

In odd characteristic, the Tate conjecture implies the Artin-Tate conjecture \cite[Theorem 6.1]{milnetate}. 
If $X$ is a K3 surface both hold \cite{tatecharles,tatepera,tatemaulik}; furthermore, 
$\# \Br(X)$ is a perfect square (see, e.g., \cite{br})
Thus,
\begin{equation}
  \label{disc}
   \disc( \NS(X_{\FF_q}) ) =  \lim_{s \rightarrow 1} \frac{  (-1)^{\rho(X) -1}   P_2(X,q^{-s})}{ q (1-q^{1-s})^{\rho(X)} } \bmod{\QQ^{\times 2}}.
\end{equation}

Usually, one computes $P_2$ by counting points in sufficiently many extensions of 
the base field. For K3 surfaces, this requires computations 
in fields of size at least $p^{10}$.
Such computations have been performed in \cite{picone,ejdoublecover1,ejdoublecover2,ejgalois,ejmodp} 
for primes $<10$. This direct approach is computationally not feasible for larger primes.
Our approach follows an idea of Kedlaya: we extract $P_2$ by
computing the Frobenius action on $p$-adic cohomology
(Monsky-Washnitzer cohomology) with sufficient precision. 
For example, for a quartic K3 surface over $\FF_p$, where $p > 41$, 
it suffices to know two significant $p$-adic digits of the coefficients of $P_2$.
This can be achieved using the Newton identities combined with Mazur inequalities \cite{mazur}.

The algorithmic implementation of this idea relies on techniques introduced in \cite{akr} and \cite{harvey}.
The approach by Abbott--Kedlaya--Roe \cite{akr} makes primes $ < 20$ computationally feasible and it was used in \cite{luijknontrivialaut}, 
but its dependence on $p$ is at least $p^{\dim(X)+1}$.
We make use of refinements of Kedlaya's algorithm, which were introduced by Harvey \cite{harvey}:
\begin{itemize}
    \item rewriting the Frobenius action on Monsky-Washnitzer cohomology in terms of sparse polynomials;
    \item preserving the sparseness throughout the reduction process of differentials in cohomology;
    \item rewriting each reduction step process as a linear map.
\end{itemize}
The time complexity is dominated by the reduction of differentials in cohomology, which
involves $O(p)$ recurrent matrix vector multiplications in $\ZZ/p^s \ZZ$. 
For a quartic K3 surface the size of the matrices is $220 \times 220$; for $p > 41$ one can take $s = 4$.
Moreover, if the K3 is nondegenerate (as in \cite{nondegenerate}), one can reduce their size to $64 \times 64$.
In practice, we had no difficulties finding  a change of coordinates for which the surface became nondegenerate.

Altogether, this reduces the polynomial dependence on $p$ in \cite{akr} to quasi-linear 
(or to $p^{1/2 + \varepsilon}$ using \cite{BGS}).
The details of the algorithm will be presented in \cite{controlledreduction}.
Our implementation is written in C++, using the libraries FLINT \cite{flint} and NTL \cite{ntl}.
The raw data of all experiments is available at {\tt www.cims.nyu.edu/$\sim$costa}.

\section{Computations and Numerical Data}

In this section, we present numerical data for five representative quartic K3 surfaces over $\QQ$ with small $\rho(\ovl{X})$.
For each surface we compute $\rho(\ovl{X})$ and 
$\rho(\ovl{X}_p)$, for all $2 <p < 2^{16}$ where $X$ 
has good reduction, using the methods introduced in Section \ref{sec:picard}.
This computation consumed around 45000 hours of CPU time of 
the Bowery and Butinah clusters at New York University.
With this data we calculate $\gamma(X,B)$ for $B < 2^{16}$, which we present as a plot.
We find two trends for $\gamma(X,B)$:
\begin{itemize}
\item
When $\rho(\overline{X}) = 1$ and $E_X = \QQ$ we have 
$$
\gamma(X,B) \sim c_X/\sqrt{B}, \quad B\ra \infty, 
$$
for some constant $c_X>0$, 
i.e., $\rho(\overline{X_p})$ jumps with probability proportional to $1/\sqrt{p}$.
\item 
When $\rho(\overline{X}) = 2$ the data leads us to believe that 
$$
\liminf_{B \rightarrow \infty} \gamma(X,B) \geq 1/2, 
$$ 
i.e., the primes at which the geometric Picard number jumps have density $\geq 1/2$.
\end{itemize}
These trends reflect which case of equation (\ref{finiteplacepicard}) we are in.

In our examples, we used the following, sufficiently generic, homogeneous polynomials:
\begin{small}
\begin{align*}
f_1 \eqdef & 2 x^2 y + 2 x y^2 + y^3 - x^2 z + x y z - y^2 z + x z^2 - 8 y z^2 + x^2 w - 9 x y w + 3 y^2 w\\
&- 10 y z w - x w^2 - 9 y w^2 + z w^2 - w^3;\\
f_2 \eqdef &-14 x^3 + x^2 y - y^3 + 2 x^2 z - 17 x y z + 22 y^2 z + x z^2 - 3 y z^2 + 2 z^3 - 2 x^2 w\\
&- 4 y^2 w - 27 x z w + y z w - 5 z^2 w - x w^2- y w^2 + 7 z w^2;\\
g_1 \eqdef& -14 x^2 - y^2 + x z + 2 y z + 2 z^2 + x w - y w - 2 z w;\\
g_2 \eqdef& -3 x^2 + 7 x y + 22 y^2 - 5 x z - z^2 - 17 x w - 27 y w + z w -4 w^2;\\
g_3 \eqdef& 2 x y + y^2 + 2 x z - y z + x w - y w + z w - w^2;\\
g_4 \eqdef& -8 x^2 + x y - y^2 - 9 y z - 9 z^2 + x w - 10 z w + 3 w^2;\\
h \eqdef &2 x^4 - 8 x^3 y - x^2 y^2 + x y^3 + y^4 + 3 x^3 z - x^2 y z + 2 x y^2 z - 10 y^3 z + x^2 z^2 - 2 x y z^2\\
&- 14 y^2 z^2 - 9 x z^3 - z^4 + x^3 w + 22 x^2 y w - 3 x y^2 w + 2 y^3 w + 7 x^2 z w + x y z w - 4 x z^2 w\\
&- 17 y z^2 w + z^3 w - 9 x^2 w^2 - x y w^2 - 5 x z w^2 - 27 y z w^2 + z^2 w^2 - y w^3 - w^4.
\end{align*}
\end{small}

We start with examples with geometric Picard number one, produced  
by forcing different lattice structures on the 
N\'eron-Severi groups on different reductions, as in \cite{picone}.

\begin{example}
\label{conicline}
Let $X$ be the smooth quartic surface over $\QQ$ defined by 
\begin{equation*}
w f_1 + p_1 z f_2 + p_2 g_1 g_2 +  p_1 p_2  h=0,
\end{equation*}
where $p_1 = 4409$ and $p_2 = 24659$.
Thus $X_{p_1}$ contains the conic $C$ defined by $w=g_1=0$, and $X_{p_2}$ contains the line $L$ defined by $w=z=0$.
Using the methods from Section \ref{sec:picard}, we find
\begin{align*}
  \rho(\overline{X}_{p_1}) = 2 & \text{ and } \disc(\NS (\overline{X}_{p_1})) = -3 \bmod{\QQ^{\times 2}};\\
  \rho(\overline{X}_{p_2}) = 2 & \text{ and } \disc(\NS (\overline{X}_{p_2})) = -1 \bmod{\QQ^{\times 2}}.
\end{align*}
Therefore, $\rho(\overline{X}) = 1$. Furthermore, $\NS (\overline{X}_{p_1})$ is generated by the hyperplane section and the conic $C$, and $\NS (\overline{X}_{p_2})$ 
is generated by the hyperplane section and the line $L$.
In this example we 
observe $\rho(\overline{X}_p) > 4$ only for $p=29$, 
where $\rho(\overline{X}_{29}) = 6$ and  $\disc(\NS (\overline{X}_{29})) = -537$.
\end{example}

\begin{example}
\label{ellipticline}
Let $X$ be the K3 surface over $\QQ$ defined by 
\begin{equation*}
p_1( w f_1 +  z f_2 ) +  p_2 (g_1 g_2 + g_3 g_4) + p_1 p_2 h=0,
\end{equation*}
with $p_1 = 18869 $ and $p_2 = 30047$.
As in the previous example, $X_{p_2}$ 
contains a line $L$. Now $X_{p_1}$ 
contains the elliptic curve $C$ defined by $g_1 = g_3 = 0$.
Nonetheless, we still have 
\begin{align*}
  \rho(\overline{X}_{p_1}) = 2 & \text{ and } \disc(\NS (\overline{X}_{p_1})) = -3 \bmod{\QQ^{\times 2}};\\
  \rho(\overline{X}_{p_2}) = 2 & \text{ and } \disc(\NS (\overline{X}_{p_2})) = -1 \bmod{\QQ^{\times 2}}.
\end{align*}
Consequently, $\rho(\overline{X}) = 1$, $\NS (\overline{X}_{p_1})$ 
is generated by the hyperplane section and the elliptic curve $C$, 
and $\NS (\overline{X}_{p_2})$ 
is generated by the hyperplane section and the line $L$.
As in the previous example, $\rho(\overline{X}_p) > 4$ for only one prime $p=7$, 
where $\rho(\overline{X}_7) = 6$ and  $\disc(\NS (\overline{X}_{7})) = -345$. 
\end{example}

\begin{figure}[h]
    \centering
    \begin{subfigure}{\graphwidthtwo}
       \includegraphics[width=\textwidth]{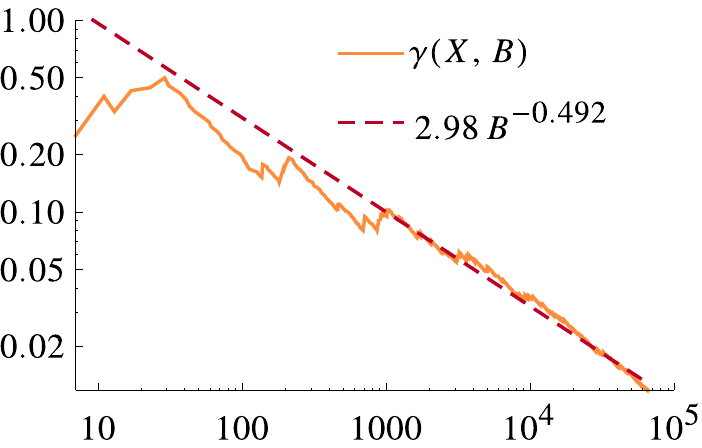} 
        \caption{Example \ref{conicline}}
    \end{subfigure}
    \begin{subfigure}{\graphwidthtwo}
        \includegraphics[width=\textwidth]{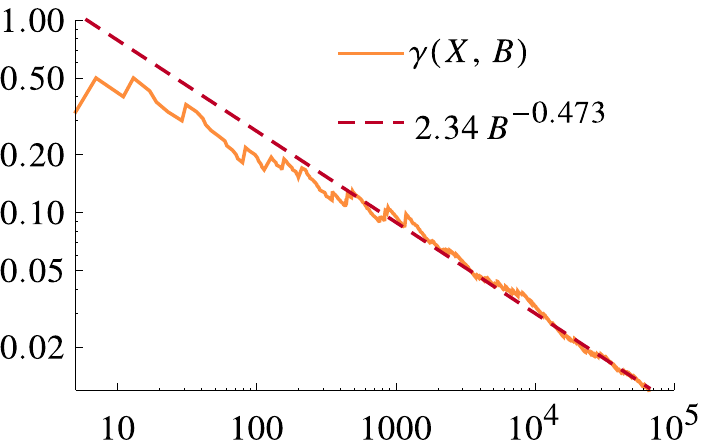}
        \caption{Example \ref{ellipticline}}
    \end{subfigure}
    \caption{Log-log plots of $\gamma$ and their least-square-fit to a power law in Examples \ref{conicline} and  \ref{ellipticline}.}
    \label{fig:picardnumberone}
\end{figure}

In both examples, $\eta(\ovl{X}) = 1$, 
$E_X = \QQ$, and $X$ does not have complex multiplication.
We present the log-log plots of $\gamma(X,B)$ for the previous examples and 
their respective least square fit to a power law in Figure \ref{fig:picardnumberone}.
We observe that
$$
\gamma(X,B) \sim \frac{c_X}{\sqrt{B}}, \quad B\ra \infty,
$$
for some constant $c_X>0$.

Next, we present examples of K3 surfaces over $\QQ$ with 
geometric Picard number two. We achieve this by forcing an additional curve on 
$X$ and by finding a prime $p$ such that $\rho(\overline{X}_p) = 2$.
\begin{example}
\label{line}
Let $X$ be the K3 surface given by 
\begin{equation*}
    w f_1 + z f_2 =0;
\end{equation*}
it contains the line $L$ defined by $w=z=0$. For $p=23$ we have $\rho(\ovl{X}_p) = 2$.
\end{example}

\begin{example}
\label{conic}
Let $X$ be the smooth quartic surface given by
\begin{equation*}
    w f_1 + g_1 g_2=0,
\end{equation*}
containing the conic $C$ defined by $w = g_1 = 0$.
For $p = 17$ we have $\rho(\overline{X}_p) = 2$.
\end{example}

\begin{example}
\label{elliptic}
Let $X$ be defined by 
\begin{equation}
g_1 g_2 + g_3 g_4=0,
\end{equation}
and containing the curve $C$ given by $g_1 = g_3 = 0$. For 
$p = 31$ we have $\rho(\ovl{X}_p) = 2$.
\end{example}

\begin{figure}[h]
    \centering
    \includegraphics[width=\graphwidthone]{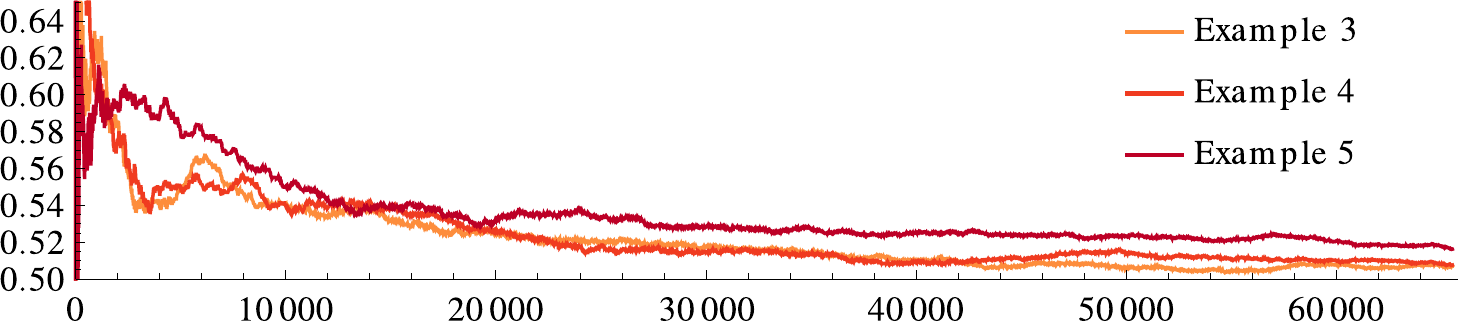}
    \caption{Plots of $\gamma(X,B)$ for the Examples \ref{line}, \ref{conic} and \ref{elliptic}.}
    \label{fig:picardnumbertwo}
\end{figure}

In Figure \ref{fig:picardnumbertwo} we present plots of $\gamma(X,B)$ for the previous examples.
These suggest that 
$$
\liminf _{B\rightarrow \infty} \gamma(X,B) \geq 1/2.
$$ 

\begin{figure}[h]
    \centering
    \includegraphics[width=\graphwidthone]{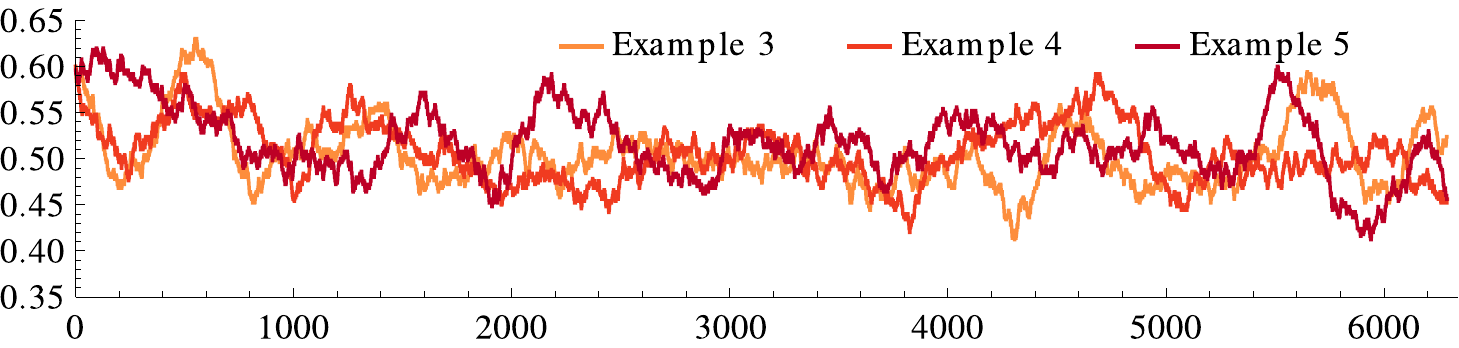}
    \caption{Plots of $\delta(X,i,250)$ for the Examples \ref{line}, \ref{conic} and \ref{elliptic}.}
    \label{fig:picardnumbertwomavg}
\end{figure}

For these examples we also inspected the local density of $\Pi_{\rm jump}(X)$. For this we resort to a moving average,
\begin{equation*}
    \delta(X,i,j) := \frac{\#\left\{ i - j/2 < l \leq i + j/2 :  \, p_l \in \Pi_{\rm jump}(X) \right\} }{j}
\end{equation*}
where $\left\{p_1,p_2,p_3,\dots\right\}$ denotes the primes, on their natural order, 
at which $X$ has good reduction.
We present $\delta(X,i,250)$ in Figure \ref{fig:picardnumbertwomavg}. 
We observe that the moving average oscillates slightly {\em above} $1/2$.

While in these examples we observed $\rho(\ovl{X}_p) > 2$ 
more frequently, the number of primes such that 
$\rho(\overline{X}_p) > 4$ is quite small. 
We present those in Table \ref{tab:geometricpicardnumbertwo}.
\begin{table}[h]
  \begin{subtable}{0.2\textwidth}
  \begin{tabular}{|r|c|}
    \multicolumn{2}{c}{Example \ref{line}}\\
\hline
$p$ & $\rho(\overline{X}_p)$\\
\hline
 3 & 6 \\
 11 & 6 \\
 13 & 6 \\
 47 & 6 \\
 53 & 6 \\
 181 & 6 \\
 239 & 6 \\
 25087 & 6 \\
\hline
\end{tabular}
\end{subtable}
\hskip 0.5cm
\begin{subtable}{0.2\textwidth}
  \begin{tabular}{|r|c|}
     \multicolumn{2}{c}{Example \ref{conic}}\\
  \hline
  $p$ & $\rho(\overline{X}_p)$\\
\hline
 3 & 10 \\
 5 & 10 \\
 11 & 6 \\
 29 & 6 \\
 83 & 6 \\
 491 & 6 \\
 2777 & 6 \\
 \phantom{0}3187 & 6 \\
\hline
\end{tabular}
\end{subtable}
\hskip 0.5cm
\begin{subtable}{0.2\textwidth}
\begin{tabular}{|r|c|}
     \multicolumn{2}{c}{Example \ref{elliptic}}\\
  \hline
  $p$ & $\rho(\overline{X}_p)$\\
  \hline
 3 & 6 \\
 17 & 6 \\
 \phantom{00}347 & 6 \\
 &  \\
&  \\
&  \\
&  \\
&  \\
  \hline
\end{tabular}
\end{subtable}

\caption{Primes $p < 2^{16}$ for which $\rho(\overline{X}_p)>4$.}
\label{tab:geometricpicardnumbertwo}
\end{table}

\bibliographystyle{alpha}
\bibliography{biblio}

\end{document}